

\font\eusm=eusm10                   


\font\eusms=eusm7                       

\font\eusmss=eusm5                      

\input amstex

\documentstyle{amsppt}
  \magnification=1000
  \hsize=7.0truein
  \vsize=9.0truein
  \hoffset -0.1truein
  \parindent=2em


\define\Ad{\text{\rm Ad}}                    




\define\alih{{(\alpha\circ\iota)\hat{\,}}}   

\define\alphah{{\hat\alpha}}                 

\define\Aut{\text{\rm Aut}}                  


\define\betat{{\widetilde\beta}}             


\define\Bof{B}                               


\define\CMp{C_{\MvN_\phi}}                   


\define\dif{\text{\it d}}                    


\define\Eh{\hat E}                           


\define\Fc{{\mathchoice                      
     {\text{\eusm F}}
     {\text{\eusm F}}
     {\text{\eusms F}}
     {\text{\eusmss F}}}}


\define\freeprodi{\mathchoice                
     {\operatornamewithlimits{\ast}
      _{\iota\in I}}
     {\raise.5ex\hbox{$\dsize\operatornamewithlimits{\ast}
      _{\sssize\iota\in I}$}\,}
     {\text{oops!}}{\text{oops!}}}

\define\gammai{{\gamma^{-1}}}                

\define\gammao{{\gamma_1}}                   

\define\gammat{{\gamma_2}}                   

\define\Gh{\hat G}                           

\define\Hil{{\mathchoice                     
     {\text{\eusm H}}
     {\text{\eusm H}}
     {\text{\eusms H}}
     {\text{\eusmss H}}}}







\define\id{\text{\it id}}                    


\define\ioh{{\hat \iota}}                    


\define\Integers{\bold Z}                    

\define\lspan{\text{\rm span}@,@,@,}         

\define\mfr{{\frak m}}                       

\define\MvN{{\Cal M}}                        

\define\MvNag{\MvN^\alpha(\{\gamma\})}       


\define\MvNpg{\MvN_\phi(\{\gamma\})}         

\define\nm#1{||#1||}                         

\define\NvN{{\Cal N}}                        

\define\NvNt{{\widetilde\NvN}}               



\define\pit{{\tilde\pi}}                     

\define\PvN{{\Cal P}}                        

\define\PvNt{{\widetilde\PvN}}               

\define\QED{$\hfill$\qed\enddemo}            

\define\QvN{{\Cal Q}}                        

\define\QvNt{{\widetilde\QvN}}               

\define\Real{{\bold R}}                      

\define\restrict{\lower .3ex                 
     \hbox{\text{$|$}}}

\define\Rp{\Real_+^*}                        

\define\RvN{{\Cal R}}                        

\define\Spec{\text{\rm Sp}}                  

\define\spw{{\widehat{\sigma^\phi}}} 

\define\supp{\text{\rm supp}}                

\define\Tcirc{\bold T}                       

\define\tr{\text{\rm tr}}                    

\define\Tr{\text{\rm Tr}}                    









\topmatter

  \title  Crossed product decompositions of a purely infinite von Neumann
          algebra with faithful, almost periodic weight \endtitle

  \rightheadtext{Crossed product decompositions}

  \author Kenneth J. Dykema \endauthor

  \thanks This work was supported by
           a National
          Science Foundation Postdoctoral Fellowship. \endthanks

  \affil Department of Mathematics \\ University of California \\
         Berkeley, CA 94720 \endaffil

  \date \enddate

  \abstract  For $\MvN$ a separable, purely infinite von Neumann algebra
             with almost periodic weight $\phi$, a decomposition of
             $\MvN$ as a crossed product of a semifinite von Neumann
             algebra by a trace--scaling action of a countable
             abelian group is given.
             Then Takasaki's continuous decomposition of the same algebra
             is related to the above discrete decomposition via
             Takesaki's notion of induced action,
             but here one induces up from a dense subgroup.
             The above results are used to give a model for the
             one--parameter trace--scaling action of $\Rp$ on
             the injective II$_\infty$ factor.
             Finally, another model of the same action, due to work
             of Aubert and explained by Jones, is described.
  \endabstract

\endtopmatter

\document
  \TagsOnRight
  \baselineskip=18pt

\noindent{\bf Introduction.}

  A crucial part of the present--day understanding of type~III factors
is their decomposition as crossed products of type~II$_\infty$ von~Neumann
algebras by groups of trace--scaling (or trace--decreasing) automorphisms.
This was accomplished by Connes in~\cite{4} and by Takesaki in~\cite{14}.
Connes defined the classification of type~III factors as type~III$_\lambda$,
$0\le\lambda\le1$ and showed
that a type III$_\lambda$ factor where $0\le\lambda<1$ is isomorphic
to the crossed product of a type II$_\infty$ von Neumann algebra, $\NvN$
by a single automorphism ({\it i.e\.} by the group $\Integers$).
When $\lambda>0$ $\NvN$ can be chosen to be a factor and the automorphism
trace--scaling, and in the case $\lambda=0$ the automorphism is ergodic
on the center of $\NvN$ and a trace can be chosen such that the
automorphism is strictly decreasing of the trace.

  Takesaki developed the theory of crossed products of
a von Neumann algebra by actions of locally compact groups,
including his duality theory and his theory of induced actions.
He thereby proved the continuous decomposition for a factor $\MvN$:
if $\MvN$ is type III$_1$ then $\MvN$ is the crossed product of
a type II$_\infty$ factor $\NvN$ by a one--parameter group of trace--scaling
automorphisms.

  Almost periodic weights
(the definition is reviewed in~\S1.3) were defined by Connes
in~\cite{4} and can be used to elucidate the structure of certain
type~III$_1$ factors.
Connes defined the
invariant $Sd$ for a full type~III$_1$ factor in terms of its
almost periodic weights in~\cite{6}, where he also showed that
there is a type~III$_1$ factor having no almost periodic weights.
However, many type~III$_1$ factors of interest have them.
For example, the injective type~III$_1$ factor,
which was shown to be unique
by Haagerup~\cite{9}, has many almost periodic weights, (cf~\S4).
Also, the free products of certain finite dimensional algebras
with respect to certain states that are not traces are known to be
type~III$_1$ factors (see~\cite{3} and~\cite{7})
and the free product states on them
are almost periodic.

  In this paper, given a purely infinite, separable
von Neumann algebra, $\MvN$, with
almost periodic weight, $\phi$, we show how $\MvN$ can
be written as a discrete crossed product of a type II$_\infty$
von Neumann algebra $\NvN$ by the action of a countable
abelian group.

  Then the continuous decomposition for $\MvN$ can be described
in terms of the above discrete decomposition via a straightforward
generalization of Takesaki's definition of induced action.  This
generalization is to induce up from a possibly dense subgroup
rather than just a closed subgroup.  One thereby obtains when
$\MvN$ is a type~III$_1$ factor a description of the core
type~II$_\infty$ factor and the one--parameter group action on it.

  Hence in the injective case a concrete description of the (unique
by Haagerup~\cite{9}) one--parameter
trace--scaling action on the (unique by Connes~\cite{5})
injective II$_\infty$ factor is obtained.
It was explained to the author by V.F.R\. Jones that there is another
description
of this action, which will be described in the last section of this paper.

  This paper has five sections.
In \S1, some preliminaries
are briefly reviewed and our notation is explained, including
\S1.1, Takesaki's definition of crossed products;
\S1.2, the Arveson--Connes spectral analysis for a compact abelian
group action;
\S1.3, Connes' compact group action on a von Neumann algebra with
almost periodic weight.
In~\S2, the discrete crossed product decomposition for a
purely infinite von Neumann algebra
with almost periodic weight is given.
In~\S3, Takesaki's continuous decomposition for a
purely infinite von Neumann algebra
with almost periodic weight is
shown to arise via an induced action from the discrete action of~\S2.
In~\S4, the previous results are used to give a concrete description
of the one--parameter trace--scaling action on the injective~II$_\infty$
factor.
In~\S5, a previously known description of this action is briefly
described.

\medpagebreak

\noindent{\bf Acknowledgements.}
  The work embodied in this paper was prompted by a desire to
understand the one--parameter trace--scaling action on the injective
II$_\infty$ factor.  This desire in turn stemmed from a discussion
with M\. Takesaki, for which I am grateful.  I would also like
to thank J\. Feldman for helpful discussions and V.F.R\. Jones
for telling me about the model in \S5.

\medpagebreak

\noindent{\bf\S1. Preliminaries and notation.}
Throughout this paper, all groups and Hilbert spaces will be separable
and all von Neumann algebras will have separable preduals.
Although these restrictions avoid measure theoretic difficulties,
it is likely that similar results hold in more general circumstances.

\noindent{\bf\S1.1. Crossed products} (Takesaki~\cite{14}).
Let $\MvN$ be a von Neumann algebra
and $\alpha:G\rightarrow\Aut(\MvN)$
a continuous action (Arveson~\cite{1}) of a locally compact group
$G$ on $\MvN$.
(Henceforth, the word ``action'' will be used to mean a continuous action.)
Suppose $\sigma$ is a normal, faithful representation of $\MvN$ as operators
on the Hilbert space $\Hil$.
Consider the von Neumann algebra
acting on $L^2(G)\otimes\Hil$ that is generated by
$$ \{\lambda_g\otimes1\mid g\in G\}\cup
   \{\pi_{\alpha,\sigma}(a)\mid a\in\MvN\},$$
where $\lambda_g\xi(h)=\xi(g^{-1}h)$ for $\xi\in L^2(G)$ and $g,h\in G$,
and where $\pi_{\alpha,\sigma}(a)\zeta(h)=\alpha_{h^{-1}}(a)\zeta(h)$
for $h\in G$ and $\zeta\in L^2(G,\Hil)$, where $L^2(G,\Hil)$ is
the space of square--integrable, measurable functions from $G$ into $\Hil$,
and is identified with $L^2(G)\otimes\Hil$.
The measure on $G$ is left Haar measure.
Takesaki proved that the choice
of representation $\sigma$ does not affect the isomorphism class, so that
one may define the {\it crossed product} von Neumann algebra, which we
will denote $\MvN\rtimes_\alpha G$, to be generated by elements denoted
$$ \{\lambda_g\otimes1\mid g\in G\}\cup
   \{\pi_{\alpha}(a)\mid a\in\MvN\}$$
in such a way that if $\sigma$ is a representation of $\MvN$ as above,
then there is a faithful, normal representation of $\MvN\rtimes_\alpha G$
given by sending $\lambda_g\otimes1$ to $\lambda_g\otimes1$ and
$\pi_\alpha(a)$ to $\pi_{\alpha,\sigma}(a)$.

  For $G$ abelian, Takesaki defined the {\it dual action} $\alphah$
of $\Gh$ on $\MvN\rtimes_\alpha G$
to be given by
$\alphah_\gamma(\lambda_g\otimes1)=\overline{\langle g,\gamma\rangle}
\lambda_g\otimes1$
and $\alphah_\gamma(\pi_\alpha(a))=\pi_\alpha(a)$.
{\it Takesaki duality} is the theorem that
$$ (\MvN\rtimes_\alpha G)\rtimes_\alphah\Gh\simeq\MvN\otimes\Bof(L^2(G)). $$

  If $\MvN$ is a purely infinite
von Neumann algebra with n.f.s\. weight $\phi$,
let $\sigma^\phi$ be the modular action of $\Real$ on $\MvN$.
Let $\Rp$ be the dual group of $\Real$ with the pairing
$\langle t,\mu\rangle=\mu^{it}$.
Then by Takesaki duality, $\MvN\simeq\MvN\otimes\Bof(L^2(\Real))
\simeq(\MvN\rtimes_{\sigma^\phi}\Real)\rtimes_\spw\Rp$.
Takesaki showed that
$\MvN_0=\MvN\rtimes_{\sigma^\phi}\Real$ is a type II$_\infty$
von Neumann algebra that admits a n.f.s\. trace $\tau$ satisfying
$\tau\circ\spw_\mu=\mu^{-1}\tau$,
and that the W$^*$--dynamical system $(\MvN_0,\Rp,\spw)$
does not depend on $\phi$ (up to isomorphism).
We will refer to this dynamical system, giving rise to
$\MvN\simeq\MvN_0\rtimes_\spw\Rp$, as
{\it Takesaki's continuous (crossed product) decomposition}
of $\MvN$.
If $\MvN$ is a type III$_1$ factor, then $\MvN_0$ is a factor
and this is the crossed product decomposition of a type III$_1$ factor
that was mentioned in the introduction.

\noindent{\bf \S1.2. The Arveson--Connes spectral analysis for a compact
abelian group action} (Arveson~\cite{1} and Connes~\cite{4}).
Let $G$ be a compact abelian group and let $\Gamma$ denote its dual group,
with pairing $\langle g,\gamma\rangle$, ($g\in G$, $\gamma\in\Gamma$).
For $\gamma\in\Gamma$ the Arveson--Connes
spectral subspace is denoted $\MvNag$.
Hence $\MvN^\alpha(\{1_\Gamma\})$ is the fixed--point subalgebra of $\MvN$
under $\alpha$, also simply denoted $\MvN^\alpha$,
and moreover
$$ \MvN^\alpha(\{\gamma_1\})\MvN^\alpha(\{\gamma_2\})
\subseteq\MvN^\alpha(\{\gamma_1\gamma_2\}) \tag1 $$
and
$$ \MvNag^*=\MvN^\alpha(\{\gamma^{-1}\}). \tag2 $$
Also, for $a\in\MvN$ let
$$ E_\gamma(a)=\int_G\langle g,\gamma\rangle\alpha_g(a)\dif g, $$
where $\dif g$ is Haar measure.

\proclaim{Lemma 1.2.1}
For $a\in\MvN$ the following are equivalent:
\roster
\item"(i)" $a\in\MvNag$\vskip1.5ex
\item"(ii)" $\alpha_g(a)=\overline{\langle g,\gamma\rangle}a$ $\forall g\in G$
\vskip1.5ex
\item"(iii)" $E_\gamma(a)=a$.
\endroster
\endproclaim
\demo{Proof}
The implication (i)$\Rightarrow$(ii) is an easy case of Connes~\cite{4},
Lemma~2.3.5 and (ii)$\Rightarrow$(i) follows directly from the
definition of $\MvNag$.
The equivalence of (ii) and (iii) is easily verified.
\QED

\proclaim{Lemma 1.2.2}
$E_\gamma$ is a normal projection of norm $1$ from $\MvN$ onto
$\MvNag$.
If $a\in\MvN$, $b\in\MvN^\alpha(\gamma_1)$ and $c\in\MvN^\alpha(\gamma_2)$ then
$$ E_\gamma(bac)=bE_{\gamma_1^{-1}\gamma\gamma_2^{-1}}(a)c. \tag3 $$
\endproclaim
\demo{Proof}
In the notation of Connes~\cite{4}~\S2.1,
$E_\gamma=\alpha(\langle\cdot,\gamma\rangle)$,
hence is a weakly continuous linear mapping from $\MvN$ into $\MvN$
and is thus normal.
It is clearly of norm~$1$.
That $E_\gamma\circ E_\gamma=E_\gamma$ is easily verified from the definition.
That $E_\gamma$ is onto $\MvNag$ follows from Lemma~1.2.1(iii).
Equation~(3) holds because $\alpha_g(bac)=b\alpha_g(a)c$.
\QED

  It makes sense to think of $E_\gamma(a)$ as the $\gamma$th term
in the Fourier series for $a$.
To see that this generalizes the usual notion of Fourier series, let
$\MvN=L^\infty(\Tcirc)$, $G=\Tcirc$ and $\alpha_z=\text{rotation by }z$.

\proclaim{Lemma 1.2.3}
The linear span of\/ $\bigcup_{\gamma\in\Gamma}\MvNag$ is a strongly dense
$*$--subalgebra of $\MvN$.
\endproclaim
\demo{Proof}
{}From~(1) and~(2) it is clear that $\bigcup\MvNag$ is a $*$--subalgebra
of $\MvN$.
Suppose $a\in\MvN$ is such that
$\Spec_\alpha(a)$ is a finite subset of $\Gamma$.
Consider the function of $G$ $f(g)=\sum_{\gamma\in\Spec_\alpha(a)}
\langle g,\gamma\rangle$.
Then $\alpha(f)a=\sum_{\gamma\in\Spec_\alpha(a)}E_\gamma(a)$,
and since the Fourier transform $\hat f$ is identically $1$
on $\Spec_\alpha(a)$,
it follows from Connes~\cite{4} Lemma~2.1.3(b) that
$a=\sum E_\gamma(a)\in\lspan\bigcup_{\gamma\in\Gamma}\MvNag$.
But the collection of $a\in\MvN$ such that $\Spec_\alpha(a)$ is finite is
by Connes~\cite{4} Lemma~2.1.4 strongly dense in $\MvN$.
\QED

\noindent{\bf \S1.3.  Connes' compact group action associated to an
almost periodic weight} (Connes~\cite{6}).
\proclaim{Definition 1.3.1}\rm (Connes~\cite{4}).
A normal, faithful, semifinite (n.f.s.)
weight $\phi$ on a von Neumann algebra $\MvN$
is {\it almost periodic} if the modular operator $\Delta_\phi$ on
$L^2(\MvN,\phi)$ is diagonalizable, {\it i.e\.} the set of eigenvectors
of $\Delta_\phi$ has dense linear span in $L^2(\MvN,\phi)$.
\endproclaim

  In this paper, when we say $\phi$ is almost periodic we will
always exclude the trivial case when the point spectrum
of $\Delta_\phi$ is $\{1\}$.

  Let $\phi$ be an almost periodic n.f.s. weight on a von
Neumann algebra $\MvN$.
Let $\Gamma\subseteq\Real_+^*$ be the subgroup generated by the
point spectrum of $\Delta_\phi$ and endowed with the discrete
topology.
Let $\ioh:\Gamma\hookrightarrow\Real_+^*$ denote the inclusion map.
Let $G$ be the compact abelian group whose dual is $\Gamma$,
with pairing denoted as in~\S1.2.
Also, consider $\Real_+^*$ to be the dual of $\Real$ under
the pairing $\langle t,\mu\rangle=\mu^{it}$ for $t\in\Real$,
$\mu\in\Real_+^*$.
Then there is a group homomorphism $\iota:\Real\rightarrow G$ determined
by $\langle\iota(t),\gamma\rangle=\langle t,\ioh(\gamma)\rangle$
$\forall\gamma\in G,\,t\in\Real$.
Note since $\ioh$ is injective that the image of $\iota$ is dense in $G$
and also that $\iota$ is injective if and only if $\Gamma$ is dense
in $\Real_+^*$.
\proclaim{Lemma 1.3.2} {\rm (Connes~\cite{6}, Proposition 1.1).}
There is an action $\alpha$ of $G$ on $\MvN$ such that
\roster
\item"(i)" $\phi\circ\alpha_g=\phi$ $\forall g\in G$
\item"(ii)" $\alpha_{\iota(t)}=\sigma_t^\phi$ $\forall t\in\Real$,
\endroster
where $\sigma^\phi$ is the modular automorphism group associated to $\phi$.
\endproclaim
\proclaim{Lemma 1.3.3}
Let $a\in\MvN$ and $\gamma\in\Gamma$.
Then $a$ belongs to the spectral subspace $\MvNag$ if and only if
$$ \phi(ba)=\ioh(\gamma)\phi(ab)\;\forall b\in\MvN. $$
\endproclaim
\demo{Proof}
This can be proved using the KMS condition exactly like in
Takesaki~\cite{13} Lemma~1.6, or see Connes~\cite{4} Lemma~3.7.5.
\QED

\noindent{\bf\S2 The discrete decomposition.}

  Let $\MvN$ be a purely infinite von Neumann algebra and $\phi$ be
a n.f.s\. almost periodic weight on $\MvN$.  Let $\alpha$ be
Connes' action of the compact group $G$ on $\MvN$ and $\Gamma=\Gh$
as in~\S1.3.
\proclaim{Definition 2.1}\rm
We have by Takesaki duality that
$\MvN\simeq(\MvN\rtimes_\alpha G)\rtimes_\alphah\Gamma$.
The dynamical system $(\MvN\rtimes_\alpha G,\Gamma,\alphah)$
is the {\it discrete decomposition} associated to $\phi$
whose {\it core} is $\MvN\rtimes_\alpha G$.
\endproclaim
This section is devoted to
elucidating the von Neumann algebra $\MvN\rtimes_\alpha G$ and the action
$\alphah$.
We will see that the core of a discrete decomposition is semifinite
and the action $\alphah$ is trace--scaling
(for the embedding of $\Gamma$ in $\Rp$).

\proclaim{Proposition 2.2}
Let $\alpha$ be an action of a compact group $G$ on a von Neumann algebra
$\MvN\subseteq\Bof(\Hil)$ and
let $\Gamma$ be the dual group of $G$.
Then the Fourier--Plancherel transform provides
an isomorphism from $\MvN\rtimes_\alpha G$ onto
the von Neumann algebra, $\QvN$,
acting on $l^2(\Gamma)\otimes\Hil$ that is generated by
$$ \{M_f\otimes1\mid f\in l^\infty(\Gamma)\}
\cup\{\lambda_\gamma\otimes a\mid\gamma\in\Gamma,a\in\MvNag\}, $$
(where $M_f$ is the multiplication operator
$(M_f\xi)(\gamma)=f(\gamma)\xi(\gamma)$ and
$(\lambda_\gamma\xi)(\gamma')=\xi(\gamma^{-1}\gamma')$
for $\xi\in l^2(\Gamma)$).
Under this isomorphism the dual automorphism $\alphah_\gamma$
on $\MvN\rtimes_\alpha G$ corresponds to the automorphism
$\Ad(\lambda_\gammai\otimes1)$ on	 $\QvN$.
\endproclaim
\demo{Proof}
By the definition of crossed product (see \S1.1) and by
Lemma~1.2.3, $\MvN\rtimes_\alpha G$
is isomorphic to the von Neumann algebra on $L^2(G)\otimes\Bof(\Hil)$
that is generated by
$$ \{\lambda_g\otimes1\mid g\in G\}\cup\{\pi_{\alpha,\id}(a)\mid
\gamma\in\Gamma,a\in\MvNag\}. \tag4 $$
But for $a\in\MvNag$, since then
$\alpha_{h^{-1}}(a)=\langle h,\gamma\rangle a$, we have
$$ \pi_{\alpha,\id}(a)=M_{\langle\cdot,\gamma\rangle}\otimes a, $$
where $\langle\cdot,\gamma\rangle$ is the function
$G\ni h\mapsto\langle h,\gamma\rangle$.
Let $\Fc:L^2(G)\rightarrow l^2(\Gamma)$ be the Fourier--Plancherel
transform given by
$(\Fc\xi)(\gamma)=\int_G\overline{\langle g,\gamma\rangle}\xi(g)\dif g$
$\forall\gamma\in\Gamma$.
We will conjugate the operators in~(4) by $\Fc\otimes1$.
Now $\Fc M_{\langle\cdot,\gamma\rangle}\Fc^{-1}=\lambda_\gamma$
and $\Fc\lambda_g\Fc^{-1}=M_{\overline{\langle g,\cdot\rangle}}$,
so $(\Fc\otimes\id)(\MvN\rtimes_\alpha G)(\Fc\otimes\id)^{-1}$
is the von Neumann
algebra on $l^2(\Gamma)\otimes\Hil$ generated by
$$ \{M_{\overline{\langle g,\cdot\rangle}}\otimes1\mid g\in G\}
\cup\{\lambda_\gamma\otimes a\mid\gamma\in\Gamma,a\in\MvNag\}. $$
But $\{\overline{\langle g,\cdot\rangle}\mid g\in G\}$ generates
$l^\infty(\Gamma)$, so conjugation by $\Fc\otimes\id$ takes
$\MvN\rtimes_\alpha G$ onto $\QvN$.

  The dual automorphism $\alphah_\gamma$ on the von Neumann algebra
generated by the set~(4) is
$\alphah_\gamma=\Ad(M_{\overline{\langle\cdot,\gamma\rangle}}
\otimes1)$.
But $\Fc M_{\overline{\langle\cdot,\gamma\rangle}}\Fc^{-1}
=\lambda_\gammai$, proving the last sentence of the proposition.
\QED

\proclaim{Remark 2.3}\rm
$\QvN\subseteq\Bof(l^2(\Gamma))\otimes\MvN$, and every element
of $x\in\Bof(l^2(\Gamma))\otimes\MvN$ can be viewed as a generalized matrix,
indexed over $\Gamma$
and having entries in $\MvN$, where the $\gamma_1,\gamma_2$ entry
$[x]_{\gamma_1,\gamma_2}$ is given by
$$ \langle[x]_{\gamma_1,\gamma_2}v,w\rangle_\Hil
=\langle x(\chi_{\{\gamma_2\}}\otimes v),\chi_{\{\gamma_1\}}\otimes w
\rangle_{l^2(\Gamma)\otimes\Hil}, \;\forall\,v,w\in\Hil, $$
where $\chi_{\{\gamma_i\}}\in l^2(\Gamma)$ is the characteristic function
of $\{\gamma_i\}$.
One can easily prove that $\QvN$ is the set of all elements
$x\in\Bof(l^2(\Gamma))\otimes\MvN$ that, when viewed as
generalized matrices, satisfy
$[x]_{\gamma_1,\gamma_2}\in\MvN^\alpha(\{\gamma_2^{-1}\gamma_1\})$
$\forall \gamma_1,\gamma_2\in\Gamma$.
One also checks that
$$ [\Ad(\lambda_\gammai\otimes1)x]_{\gamma_1,\gamma_2}
=[(\lambda_\gammai\otimes1)x(\lambda_\gamma\otimes1)]_{\gamma_1,\gamma_2}
=[x]_{\gamma\gamma_1,\gamma\gamma_2}. $$
\endproclaim

\proclaim{Proposition 2.4}
Let $\MvN\subseteq\Bof(\Hil)$ be a von Neumann algebra and $\phi$
a n.f.s\. almost periodic weight on $\MvN$.
Let $G$, $\alpha$ and $\Gamma$ be as in~\S1.3, and let $\QvN$
be as in Proposition~2.2.
Then there is a n.f.s\. trace $\tr_\QvN$ on $\QvN$ such that
$$ \tr_\QvN\circ\alphah_\gamma=\ioh(\gamma)^{-1}\tr_\QvN. $$.
\endproclaim
\demo{Proof}
Since $\QvN\subseteq\Bof(l^2(\Gamma))\otimes\MvN$ let $\tr_\QvN$
be $\Tr(M_\ioh\cdot)\otimes\phi$ restricted to $\QvN$,
where $\Tr$ is the n.f.s\. trace on $\Bof(l^2(\Gamma))$
and the density matrix $M_\ioh$ is the multiplication operator
associated to the unbounded positive function $\ioh$ on $\Gamma$.
Then $\tr_\QvN$ is normal and semifinite because $\Tr$ and $\phi$
are n.f.s\. and $M_\ioh$ is affiliated to $\QvN$.
To see that $\tr_\QvN$ is a trace, let
$\mfr=\{a\in\MvN\mid\phi(a^*a)<+\infty\}$.
Then there is a net of projections
$e_i\in\MvN^\alpha$, ($i\in I$), increasing to one,
that together with Lemma~1.2.3 shows that
$\lspan\bigcup_{\gamma\in\Gamma}(\MvNag\cap\mfr)$ is a strongly
dense $*$--subalgebra of $\MvN$.
Hence it will suffice to show that $\tr_\QvN(xy)=\tr_\QvN(yx)$
whenever $x$ and $y$ are finite products of elements in
$$ \{M_f\mid f\in l^\infty(\Gamma)\text{ having finite support}\}
\cup\{\lambda_\gamma\otimes a\mid\gamma\in\Gamma,\,a\in\MvNag\cap\mfr\}. $$
But such a finite product is equal to
$(M_f\otimes1)(\lambda_\gamma\otimes a)$ for some $f\in l^\infty(\Gamma)$
having finite support, some $\gamma\in\Gamma$ and $a\in\MvNag$.
But this in turn is the sum of operators of the form
$(\chi_{\gamma'}\otimes1)(\lambda_\gamma\otimes a)$, which when viewed
as a generalized matrix as in Remark~2.3 has all entries equal to zero
except the $\gamma',\gamma'\gamma^{-1}$ entry.
Hence it suffices to show that $\tr_\QvN(xy)=\tr_\QvN(yx)$ when
$[x]_{\gamma_1,\gamma_2}=a$ and all other entries of $x$ are zero
and $[y]_{\gamma_3,\gamma_4}=b$ and all other entries of $y$ are zero.
Both $\tr_\QvN(xy)$ and $\tr_\QvN(yx)$ are zero unless $\gamma_3=\gamma_2$
and $\gamma_4=\gamma_1$, so assume this is the case.
Then $a\in\MvN^\alpha(\{\gamma_2^{-1}\gamma_1\})$ and
$b\in\MvN^\alpha(\{\gamma_1^{-1}\gamma_2\})$, so by Lemma~1.3.3,
$$ \tr_\QvN(yx)=\ioh(\gamma_2)\phi(ba)=\ioh(\gamma_1)\phi(ab)=\tr_\QvN(xy). $$
\QED

  Let us now take a closer look at spectral subspaces $\MvNag$.
Taking $\MvN$, $\phi$, {\it etc\.} as in the previous proposition,
by Lemma 1.3.3 $\MvN^\alpha$ equals the centralizer, $\MvN_\phi$,
of $\phi$.
Thus for $\gamma\in\Gamma$, we will also denote the spectral subspace
$\MvNag$ by $\MvNpg$.
Let $Z_\phi=Z(\MvN_\phi)$ denote the center of the centralizer of $\phi$
and $PZ_\phi$ the set of (self--adjoint) projections in $Z_\phi$.
In general, for $a$ an element in a von Neumann algebra $\NvN$,
$C_\NvN(a)$ will denote the central carrier of $a$ in $\NvN$.
\proclaim{Definition 2.5}\rm
For $\gamma\in\Gamma$ and $p\in PZ_\phi$ let
$$ S_\gamma(p)=\bigvee\{\CMp(apa^*)\mid a\in\MvNpg\}\in PZ_\phi. $$
\endproclaim

\proclaim{Lemma 2.6}
\roster
\item"(i)" $\forall\gamma\in\Gamma$ there is a partial isometry
$v_\gamma\in\MvNpg$ such that
$$ S_\gamma(p)=\CMp(v_\gamma pv_\gamma^*)\;\forall p\in PZ_\phi; $$
\item"(ii)" one may arrange that $v_\gammai=v_\gamma^*$.
\item"(iii)" Let $\gamma\in\Gamma$;
then $p\le p'\Rightarrow S_\gamma(p)\le S_\gamma(p')$
$\forall p,p'\in PZ_\phi$;
\item"(iv)" setting $p_\gamma=S_\gammai(1)$ we have
$S_\gamma(p)=0$ if $p\perp p_\gamma$ and
$$ S_\gammai S_\gamma(p)=pp_\gamma\;\forall p\in PZ_\phi; $$
\item"(v)" if $p_1,p_2\in p_\gamma PZ_\phi$ then
$S_\gamma(p_1)=S_\gamma(p_2)$ implies $p_1=p_2$;
\item"(vi)" $S_\gamma(p_1p_2)=S_\gamma(p_1)S_\gamma(p_2)$
$\forall p_1,p_2\in PZ_\phi$;
\item"(vii)" if $p_1,p_2\in PZ_\phi$ and $p_1\perp p_2$ then
$S_\gamma(p_1+p_2)=S_\gamma(p_1)+S_\gamma(p_2)$.
\item"(viii)" $S_\gamma$ is normal in the sense that if $(p_n)_{n=1}^\infty$
is an increasing family of projections in $Z_\phi$ then
$$ S_\gamma(\bigvee_{n=1}^\infty p_n)=\bigvee_{n=1}^\infty S_\gamma(p_n); $$
\endroster
\endproclaim
\demo{Proof}
For~(i), let $(v_i)_{i\in I}$ be a family of partial isometries in $\MvNpg$
that is maximal with respect to the property that
$\CMp(v_i^*v_i)\perp\CMp(v_{i'}^*v_{i'})$ and
$\CMp(v_iv_i^*)\perp\CMp(v_{i'}v_{i'}^*)$ whenever $i,i'\in I$, $i\neq i'$.
Let $v_\gamma=\sum_{i\in I}v_i$.
For $p\in PZ_\phi$, let $q=\CMp(v_\gamma pv_\gamma^*)$.
Clearly $q\le S_\gamma(p)$.
Suppose for contradiction that $q\neq S_\gamma(p)$.
Then $\exists a\in\MvNpg$ such that $(1-q)ap\neq0$.
Taking the polar decomposition, let $w$ be the polar part of $(1-q)ap\neq0$.
Then $w\in\MvNpg$ is a partial isometry, $w^*w\le p$ and $ww^*\le(1-q)$.
We will show that $\CMp(ww^*)\perp\CMp(v_iv_i^*)$ and
$\CMp(w^*w)\perp\CMp(v_i^*v_i)$ $\forall i\in I$, which will contradict
the maximality of $(v_i)_{i\in I}$.
If for some $i\in I$, $\CMp(ww^*)\CMp(v_iv_i^*)\neq0$ then
$\exists a\in\MvN_\phi$ such that $w^*av_i\neq0$.
But $w^*av_i\in\MvN_\phi$ so $0\neq pw^*(1-q)av_i=w^*a(1-q)v_ip$
and $(1-q)v_ip\neq0$, contradicting the choice of $q$.
Similarly, if if $\CMp(w^*w)\CMp(v_i^*v_i)\neq0$ then
$\exists a\in\MvN_\phi$ such that $wav_i^*\neq0$, so
$0\neq(1-q)wpav_i^*=wapv_i^*(1-q)$, a contradiction.

  For~(ii), just note from the proof of~(i) that any maximal family
$(v_i)_{i\in I}$ will do, and use that
$\MvN_\phi(\{\gamma^{-1}\})=\MvNpg^*$.

  Part~(iii) is clear from the definition of $S_\gamma$.  For~(iv),
we have from~(1) and~(ii) that $p_{\gamma}=\CMp(v_\gamma^*v_\gamma)$.
If $p\in PZ_\phi$ and $p\perp p_\gamma$, then $p\perp v_\gamma^*v_\gamma$,
so $S_\gamma(p)=\CMp(v_\gamma pv_\gamma^*)=0$.
For any $p\in PZ_\phi$, $S_\gamma(p)\ge v_\gamma pv_\gamma^*$ so
$S_\gammai S_\gamma(p)=\CMp(v_\gamma^* S_\gamma(p)v_\gamma)
\ge\CMp(pv_\gamma^*v_\gamma)=pp_\gamma$.
But also $S_\gammai S_\gamma(p)=\CMp(v_\gamma^* S_\gamma(p)v_\gamma)
=\CMp(v_\gamma^*\CMp(v_\gamma pv_\gamma^*)v_\gamma)
\le\CMp(v_\gamma^*v_\gamma pv_\gamma^*v_\gamma)=pp_\gamma$,
proving~(iv).  Part~(v) is clear from~(iv).

  For~(vi), $S_\gamma(p_1p_2)\le S_\gamma(p_1)S_\gamma(p_2)$ from~(iii).
So using~(iv), $p_\gamma p_1p_2=S_\gammai S_\gamma(p_1p_2)
\le S_\gammai(S_\gamma(p_1)S_\gamma(p_2))
\le(S_\gammai S_\gamma(p_1))(S_\gammai S_\gamma(p_2))=p\gamma p_1p_2$.
Hence the inequalities in the previous sentence are equalities,
so by~(v), $S_\gamma(p_1p_2)=S_\gamma(p_1)S_\gamma(p_2)$.

  For~(vii), we have from~(vi) and~(i) that
$\CMp(v_\gamma p_1v_\gamma^*)\perp\CMp(v_\gamma p_2v_\gamma^*)$,
hence $S_\gamma(p_1+p_2)=\CMp(v_\gamma p_1v_\gamma^*+v_\gamma p_2v_\gamma^*)
=\CMp(v_\gamma p_1v_\gamma^*)+\CMp(v_\gamma p_2v_\gamma^*)
S_\gamma(p_1)+S_\gamma(p_2)$.

  Part~(viii) holds because taking central carriers respects $\bigvee$
(see Kadison and Ringrose~\cite{10}, 5.5.3).
\QED

\proclaim{Lemma 2.7}
\roster
\item"(i)" $S_\gamma$ extends to a normal $*$--homomorphism,
also denoted $S_\gamma$, from $Z_\phi$ into $Z_\phi$;
the kernel of $S_\gamma$ is $(1-p_\gamma)Z_\phi$ and the restriction of
$S_\gamma$ is an isomorphism from $p_\gamma Z_\phi$ onto
$p_\gammai Z_\phi$.
\item"(ii)" We have $p_{1_\Gamma}=1$ and $S_{1_\Gamma}=\id$.
\item"(iii)" Let $\gamma_1,\gamma_2\in\Gamma$ and
$q=S_{\gamma_2^{-1}}(p_{\gamma_1})$.
Then $q\le p_{\gamma_1\gamma_2}$ and
$$ S_{\gamma_1}S_{\gamma_2}\restrict_{qZ_\phi}
=S_{\gamma_1\gamma_2}\restrict_{qZ_\phi}. $$
\endroster
\endproclaim
\demo{Proof}
For~(i), the properties proved in Lemma~2.6 show that $S_\gamma$
extends to a homomorphism from $\lspan PZ_\phi$ ({\it i.e\.} the $L^\infty$
functions taking only finitely many values) to $\lspan PZ_\phi$, with
kernel $\lspan(1-p_\gamma)PZ_\phi$ and that is 1--1 on
$\lspan p_\gamma PZ_\phi$.  Extending to all of $PZ_\phi$ is now standard
measure theory, making use of part~(viii) of Lemma~2.6.
Part~(ii) is clear.

  For~(iii), let us show that
$$ S_{\gamma_1}S_{\gamma_2}(p)=S_{\gamma_1\gamma_2}(qp)\;
\forall p\in PZ_\phi. \tag5 $$
Denoting by $\supp(b)$ the support in $\MvN_\phi$ of a self--adjoint element
$b\in\MvN_\phi$, we have
$$ \spreadlines{2.5ex} \align
S_{\gamma_1}S_{\gamma_2}(p)
&=\bigvee_{a_1\in\MvN_\phi(\{\gamma_1\})}\supp(a_1
(\bigvee_{a_2\in\MvN_\phi(\{\gamma_2\})}\supp(a_2pa_2^*))a_1^*) \\
&\le\bigvee\Sb a_1\in\MvN_\phi(\{\gamma_1\}) \\
a_2\in\MvN_\phi(\{\gamma_2\}) \endSb a_1a_2pa_2^*a_1^* \\
&\le\bigvee_{a\in\MvN_\phi(\{\gamma_1\gamma_2\})}apa^*
=S_{\gamma_1\gamma_2}(p).
\endalign $$
Hence
$$ S_\gammao S_\gammat(p)=S_\gammao(p_\gammao S_\gammat(p))
=S_\gammao(p_\gammao p_{\gamma_2^{-1}}S_\gammat(p))
=S_\gammao S_\gammat(qp)\le S_{\gammao\gammat}(qp). $$
But $qp=S_{\gamma_2^{-1}}(p_\gammao S_\gammat(p))$ so
$$ S_{\gammao\gammat}(qp)\le S_\gammao(p_\gammao S_\gammat(p))
=S_\gammao S_\gammat(p), $$
which proves~(5).

  Now it only remains to show that $q\le p_{\gamma_1\gamma_2}$.
For this it will suffice by part~(i) to show that $p\in PZ_\phi$,
$p\le q$ and $S_{\gammao\gammat}(p)=0$ implies $p=0$.
But then by~(5) $S_\gammao S_\gammat(p)=0$, so by part~(i)
$S_\gammat(p)\perp p_\gammao$.
Hence (because $p\le q\le p_\gammat$) by Lemma~2.6iv,
$p=S_{\gamma_2^{-1}}S_\gammat(p)\perp S_{\gamma_2^{-1}}(p_\gammao)=q$,
so $p=0$.
\QED

\proclaim{Remark 2.8}\rm
{}From the above lemmas we see by standard arguments that
$Z_\phi=L^\infty(X,\mu)$ where $X$
is a standard Borel space with measure $\mu$ and for each $\gamma\in\Gamma$,
$p_\gamma$ is the characteristic function of $A_\gamma\subseteq X$ and
$S_\gamma(f)=f\circ T_\gammai$ where $T_\gammai:A_\gammai\rightarrow A_\gamma$
is a bijective measurable mapping.
Moreover, $T_{\gamma^{-1}}=(T_\gamma)^{-1}$ and
$T_\gammat\circ T_\gammao\restrict_B=T_{\gamma_2\gamma_1}\restrict_B$
where $B=T_{\gamma_1^{-1}}(A_{\gamma_1^{-1}}\cap A_\gammat)
\subseteq A_{\gammao \gammat}$.

  Thus we can describe $\Gamma\ni\gamma\mapsto T_\gamma$ as
a {\it partial action} of $\Gamma$ on $X$.
Let us define a measurable subset $F\subseteq X$ to be {\it invariant}
under the partial action $T$ if $T_\gamma(F\cap A_\gamma)=F\cap A_\gammai$
$\forall\gamma\in\Gamma$, and the partial action $T$ to be {\it ergodic}
if every invariant subset is null or conull.
An equivalent condition for the ergodicity of $T$ is in terms of
measurable equivalence relations of Feldman and Moore~\cite{8}:
the union of the graphs of $T_\gamma$ as $\gamma$ ranges over
all of $\Gamma$ is a measurable equivalence relation.
Then $T$ is ergodic if and only if the equivalence relation is ergodic.
\endproclaim

  The following proposition says that $T$ is ergodic
if and only if $\MvN$ is a factor;
but in order to avoid using the techniques of measure theory,
the condition of ergodicity is expressed in terms of the maps
$S_\gamma$.

\proclaim{Proposition 2.9}
Let $x\in\MvN$.
Then $x$ is in the center of $\MvN$ if and only if
$x$ is in the center of $\MvN_\phi$ and
$$ S_\gamma(x)=xp_\gammai\;\forall\gamma\in\Gamma. \tag6 $$
\endproclaim
\demo{Proof}
Since both $Z(\MvN)$ and the set of $x\in\MvN$ satisfying~(6)
are von Neumann subalgebras of $\MvN$, we may assume $x$ is a
projection in $\MvN$.
To show necessity, suppose $p\in PZ(\MvN)$.
Then $\phi(pb)=\phi(bp)$ $\forall b\in\MvN$ so $p\in\MvN_\phi$,
thus $p\in PZ_\phi$.
Moreover for all $\gamma\in\Gamma$ we have
$$ S_\gamma(p)=\bigvee\{\CMp(apa^*)\mid a\in\MvNpg\}=\bigvee\{p\CMp(aa^*)
\mid a\in\MvNpg\}=pp_\gammai. $$

To show sufficiency, let $x=p\in PZ_\phi$ such that~(6) holds.
To show $p\in Z(\MvN)$ it will by Lemma~1.2.3 suffice to show that
$p$ commutes with $a$ for every $\gamma\in\Gamma$ and $a\in\MvNpg$.
But if not then $(1-p)ap\neq0$, so $(1-p)\CMp(apa^*)\neq0$,
so $(1-p)S_\gamma(p)\neq0$, contradicting $S_\gamma(p)=pp_\gammai$.
\QED

\proclaim{Corollary 2.10}
$\MvN$ is a factor if and only if the partial action $T$ is ergodic.
\endproclaim

\proclaim{Proposition 2.11}
Let $\MvN$, $\phi$, $G$, $\alpha$, $\Gamma$ and $\QvN\simeq\MvN\rtimes_\alpha
G$
be as in Proposition~2.4, and regard elements of $\QvN$ as generalized matrices
as described in Remark~2.3.
Let $x\in \QvN$.
Then $x$ is in the center of $\QvN$ if and only if
$$\align
\text{(1) }\;\;& [x]_{\gammao,\gammat}=0\text{ if }\gammao\neq\gammat; \\
\text{(2) }\;\;& [x]_{\gamma,\gamma}\in Z_\phi\;\forall \gamma\in\Gamma; \\
\text{and (3) }\;\;& S_\gammao([x]_{\gamma,\gamma}p_\gammao)
=[x]_{\gamma\gammao,\gamma\gammao}p_{\gamma_1^{-1}}\;
\forall \gamma,\gammao\in\Gamma.
\endalign $$
\endproclaim
\demo{Proof}
The set of $x\in\QvN$ satisfying (1)--(3) is a von Neumann subalgebra
of $\QvN$, so we may assume $x=q$ is a projection in $\QvN$.
To prove necessity, suppose $q\in Z(\QvN)$.
Then since every bounded, diagonal operator having entries in $\MvN_\phi$ is
in $\QvN$, clearly~(1) and ~(2) hold.
Suppose for contradiction that~(3) fails for
some $\gamma_1,\gamma_2\in\Gamma$.
We may assume that
$(1-[q]_{\gamma\gammao,\gamma\gammao})S_\gammao([q]_{\gamma,\gamma}p_\gammao)
\neq0$,
(otherwise take $S_{\gamma_1^{-1}}$ of both sides).
Hence $\exists a\in\MvN_\phi(\{\gammao\})$ such that
$a[q]_{\gamma,\gamma}=a\neq0$ and
$aa^*=a[q]_{\gamma,\gamma}a^*\perp[q]_{\gamma\gammao,\gamma\gammao}$.
Let $y\in\QvN$ have $(\gamma\gammao,\gamma)$ entry equal to $a$ and
all other entries equal to $0$.
It follows that $[yq]_{\gamma\gammao,\gamma}=a[q]_{\gamma,\gamma}\neq0$
and $[qy]_{\gamma\gammao,\gamma}=[q]_{\gamma\gammao,\gamma\gammao}a=0$,
contradicting that $q\in Z(\QvN)$.

  For sufficiency, suppose $q\in P\QvN$ satisfies (1)--(3).
To show $q\in Z(\QvN)$ it will suffice to show that $yq=qy$
whenever $y\in\QvN$ has only one nonzero entry, {\it i.e\.} to show
$$ a[q]_{\gamma,\gamma}=[q]_{\gamma\gammao,\gamma\gammao}a
\;\forall\gamma,\gammao\in\Gamma\;
\forall a\in\MvN_\phi(\{\gamma_1\}). \tag7 $$
Since $a^*a\in\MvN_\phi$, taking the polar decomposition of $a$,
one sees that it suffices to show~(7) for $a=v\in\MvN_\phi(\{\gammao\})$
a partial isometry.
But $vv^*\le p_{\gamma_1^{-1}}$, so using~(3) it suffices to show
$v[q]_{\gamma,\gamma}v^*=S_\gammao([q]_{\gamma,\gamma}p_\gamma)vv^*$,
hence since $v^*v\le p_\gamma$ to show
$$ vrv^*=S_\gammao(r)vv^*\;\forall r\in p_\gammao PZ_\phi\;
\forall\gammao\in\Gamma\;
\forall\text{ partial isometries }v\in\MvN_\phi(\{\gammao\}). \tag8 $$
But $vrv^*\le S_\gammao(r)$ by definition, so~$\le$ holds in~(8).
Since $r=S_{\gamma_1^{-1}}S_\gammao(r)$ we may apply~$\le$ of~(8)
to get
$vrv^*=v(v^*vS_{\gamma_1^{-1}}(S_\gammao(r)))v^*\ge
v(v^*S_\gammao(r)v)v^*=S_\gammao(r)vv^*$.
\QED

\proclaim{Proposition 2.12}
$\QvN$ is a factor if and only if $\MvN_\phi$ is a factor, and in that
case $\Gamma$ equals the point spectrum of $\Delta_\phi$
and $\QvN\simeq\MvN_\phi\otimes\Bof(\Hil)$,
where $\Hil$ is separable, infinite dimensional Hilbert space.
\endproclaim
\demo{Proof}
Suppose that $\MvN_\phi$ is not a factor
and let $p$ be a nontrivial projection
in the center of $\MvN_\phi$ and let $q\in\QvN$, satisfying (1)--(3)
of the proposition, be such that
$[q]_{\gamma,\gamma}=S_\gamma(p)$.
Then by Proposition~2.11, $q$ is in the center of $\QvN$ and $\QvN$
is not a factor.

  Suppose that $\MvN_\phi$ is a factor.
Then for each $\gamma\in\Gamma$, $p_\gamma=0$ or~$1$,
and $p_\gamma=0$ if and only if $\MvNpg=\{0\}$
if and only if $\gamma$ is not in the point spectrum of $\Delta_\phi$.
By Lemma~2.7iii, $\Gamma'=\{\gamma\in\Gamma\mid p_\gamma=1\}$ is a subgroup of
$\Gamma$.
But $\Gamma'$ equals the point spectrum of $\Delta_\phi$ and $\Gamma$ was taken
to be the group generated by $\Gamma'$, so $\Gamma'=\Gamma$.
Choose any $\gamma_0\in\Gamma$ and let $q\in\QvN$ be such that
$[q]_{\gamma_0,\gamma_0}=1$ and all other entries of $q$ are zero.
Then since $p_\gamma=1$ $\forall\gamma\in\Gamma$,
it follows that the central carrier of $q$ in $\QvN$ equals the identity
of $\QvN$.
But $q\QvN q\simeq\MvN_\phi$, so $\QvN$ is a factor.
Looking to Proposition~2.4 and noting that $\Gamma$ is an infinite set,
we see that $\QvN$ is a type II$_\infty$ factor and
$\QvN\simeq q\QvN q\otimes\Bof(\Hil)$.
\QED

\noindent{\bf \S3. The continuous decomposition}

  In this section, for $\MvN$ a purely infinite
von Neumann algebra with n.f.s\.
almost periodic weight $\phi$, Takesaki's continuous
decomposition for $\MvN$ (see~\S1.1) is related
to the discrete decomposition obtained in~\S2.
We shall see that Takesaki's continuous decomposition can be viewed
as the induced representation of the discrete decomposition, if
one broadens Takesaki's definition of induced representation
to include actions induced up from dense subgroups.
What follows is Takesaki's definition~\cite{14} of induced action,
except that he required $H$ to be a closed subgroup.
\proclaim{Definition 3.1}\rm
Let $K$ be a locally compact abelian group, $H\subseteq K$ a subgroup.
Suppose $\alpha$ is an action of $H$ on a von Neumann algebra $\NvN$.
Consider the von Neumann algebra $L^\infty(K)\otimes\NvN\simeq
L^\infty(K,\NvN)$, equal to the set of bounded measurable
functions from $K$ into $\NvN$, and the action $\beta$ of $H$
on $\L^\infty(K,\MvN)$ given by $(\beta_hf)(k)=\alpha_h(f(kh))$
for $h\in H$, $f\in\L^\infty(K,\NvN)$ and $k\in K$.
Let $\MvN$ be the fixed--point subalgebra of $L^\infty(K,\NvN)$
under the action $\beta$ of $H$.
The action, $\tau$, of $K$ on $L^\infty(K,\NvN)$ given by
$(\tau_kf)(k')=f(k^{-1}k')$ leaves $\MvN$ globally invariant.
The action {\it induced up} to $K$ from the action $\alpha$
of $H$ on $\NvN$ is the action $\tau$ of $K$ restricted to $\MvN$.
\endproclaim

  Let $G$ be a compact abelian group with dual group $\Gamma$,
$\alpha$ an action of $G$ on
a von Neumann algebra $\MvN$, $E$ a locally compact abelian group,
$\iota:E\rightarrow G$ a continuous homomorphism whose image is
dense in $G$.
Then $\alpha\circ\iota$ is an action of $E$ on $\MvN$.
Let $\ioh:\Gamma\rightarrow\Eh$ be the injective homomorphism
given by $\langle t,\ioh(\gamma)\rangle=\langle\iota(t),\gamma\rangle$
$\forall t\in E,\,\gamma\in\Gamma$.
Let
$$\align
\PvN&=\MvN\rtimes_\alpha G \\
\PvNt&=\MvN\rtimes_{\alpha\circ\iota}E. \\
\endalign $$
\proclaim{Proposition 3.2}
Let $\pit$ be a normal, faithful representation of $\MvN$ on the
Hilbert space $\Hil$.
Then the Fourier--Plancherel transform provides an isomorphism from
$\PvNt$ onto the von Neumann algebra,
$\QvNt$, acting on $L^2(\Eh)\otimes\Hil$, that is generated by
$$ \{M_F\otimes1\mid F\in L^\infty(\Eh)\}
\cup\{\lambda_{\ioh(\gamma)}\otimes\pit(a)\mid\gamma\in\Gamma,a\in\MvNag\}, $$
(where $M_F$ is the multiplication operator on $L^2(\Eh)$).
Under this isomorphism, the dual automorphism
$\alih_\mu$ on $\PvNt$ for $\mu\in\Eh$
corresponds to $\Ad(\lambda_{\mu^{-1}}\otimes1)$ on $\QvNt$.
\endproclaim
\demo{Proof}
This proposition is proved just like Proposition~2.2, once one notes
that if $a\in\MvNag$ and $t\in E$ then $(\alpha\circ\iota)_{t^{-1}}(a)
=\langle t,\ioh(\gamma)\rangle a$, so that
$\pi_{\alpha\circ\iota,\pit}(a)
=M_{\langle\cdot,\ioh(\gamma)\rangle}\otimes\pit(a)$.
\QED

\proclaim{Proposition 3.3}
In the situation of Proposition 3.2,
let the action of $\Eh$ that is induced up from the action $\alphah$
of $\Gamma$ on $\PvN$ be denoted the action $\tau$ of $\Eh$ on $\RvN$.
Then there is an isomorphism from $\PvNt$ onto $\RvN$ that intertwines
the dual action $\alih$ of $\Eh$ on $\PvNt$ with the action $\tau$.
\endproclaim
\demo{Proof}
Let $\pi_{\alpha,\id}$ be the representation of $\MvN$ on
$L^2(G)\otimes\Hil$ used in the proof of Proposition~2.2
and let $\Fc:L^2(G)\rightarrow l^2(\Gamma)$ be the Fourier--Plancherel
transform.
Then $\pit=(\Fc\otimes1)\pi_{\alpha,\id}(\Fc\otimes1)^{-1}$
is a faithful, normal representation of $\MvN$ on $l^2(\Gamma)\otimes\Hil$
and $\pit(a)=\lambda_\gamma\otimes a$ if $\gamma\in\Gamma$ and
$a\in\MvNag$, so by Proposition~3.2, $\PvNt$ is isomorphic to
the von Neumann algebra acting on $L^2(\Eh)\otimes l^2(\Gamma)\otimes\Hil$
that is generated by
$$ \{M_F\otimes1\otimes1\mid F\in L^\infty(\Eh)\}
\cup\{\lambda_{\ioh(\gamma)}\otimes\lambda_\gamma\otimes a
\mid\gamma\in\Gamma,\,a\in\MvNag\} $$
and the dual automorphism
$(\alpha\circ\iota)\hat{\,}_\mu$
for $\mu\in\Eh$ is given by $\Ad(\lambda_{\mu^{-1}}
\otimes1\otimes1)$.

  By Proposition~2.2, $L^\infty(\Eh,\PvN)$ is isomorphic to the
von Neumann algebra acting on $L^2(\Eh)\otimes l^2(\Gamma)\otimes\Hil$
that is generated by
$$ \{M_{\langle t,\cdot\rangle}\otimes1\otimes1\mid t\in E\}
\cup\{1\otimes M_{\langle g,\cdot\rangle}\otimes1\mid g\in G\}
\cup\{1\otimes\lambda_\gamma\otimes a\mid\gamma\in\Gamma,\,a\in\MvNag\},
$$
the automorphism $\beta_\gamma$ for $\gamma\in\Gamma$ is given by
$\Ad(\lambda_{\ioh(\gamma)^{-1}}\otimes\lambda_\gammai\otimes1)$
and the automorphism $\tau_\mu$ for $\mu\in\Eh$ is given by
$\Ad(\lambda_\mu\otimes1\otimes1)$.
Hence
$$
\RvN
=\left(\{M_{\langle t,\cdot\rangle}\otimes M_{\langle g,\cdot\rangle}\otimes1
\mid t\in E,\,g\in G\}
\cup\{1\otimes\lambda_\gamma\otimes a
\mid\gamma\in\Gamma,\,a\in\MvNag\}\right)''
\cap\{\lambda_{\ioh(\gamma)^{-1}}\otimes\lambda_{\gamma^{-1}}\otimes1
\mid\gamma\in\Gamma\}'.
$$
Let $U$ be the unitary on $L^2(\Eh\times\Gamma)$ given by
$(U\xi)(\mu,\gamma)=\xi(\ioh(\gamma)\mu^{-1},\gamma)$
for $\mu\in\Eh$, $\gamma\in\Gamma$,
so $U^*=U$.
Then
$$ \align
U(M_{\langle t,\cdot\rangle}\otimes M_{\langle g,\cdot\rangle})U
&=M_{\overline{\langle t,\cdot\rangle}}\otimes
M_{\langle\ioh(t)g,\cdot\rangle}, \\
U(1\otimes\lambda_\gamma)U&=\lambda_{\ioh(\gamma)}\otimes\lambda_\gamma, \\
U(\lambda_{\ioh(\gamma)^{-1}}\otimes\lambda_{\gamma^{-1}})U
&=1\otimes\lambda_\gammai, \\
U(\lambda_\mu\otimes1)U&=\lambda_{\mu^{-1}}\otimes1
\endalign $$
and
$$ \align U\RvN U
=&\left(\{M_F\otimes M_f\otimes1
\mid F\in L^\infty(\Eh),\,f\in l^\infty(\Gamma)\}
\cup\{\lambda_{\ioh(\gamma)}\otimes\lambda_\gamma\otimes a
\mid\gamma\in\Gamma,\,a\in\MvNag\}\right)'' \\
&\cap\{1\otimes\lambda_\gammai\otimes1\mid\gamma\in\Gamma\}' \\
=&\left(\{M_f\otimes1\otimes1\mid F\in L^\infty(\Eh)\}
\cup\{\lambda_{\ioh(\gamma)}\otimes\lambda_\gamma\otimes a
\mid\gamma\in\Gamma,\,a\in\MvNag\}\right)''.
\endalign $$
Moreover, the automorphism $\tau_\mu$ of $\RvN$ corresponds to the
automorphism $\Ad(\lambda_{\mu^{-1}}\otimes1\otimes1)$ of $U\RvN U$.
This is precisely the picture for $\PvNt$ and $(\alpha\circ\iota)\hat{\,}_\mu$
obtained above.
\QED

\proclaim{Corollary 3.4}
Let $\MvN$ be a purely infinite von Neumann algebra having separable
predual and $\phi$ a n.f.s\. almost periodic weight on $\MvN$.
Let $\Gamma$ be the subgroup of $\Rp$ generated by the point
spectrum of $\Delta_\phi$ and let $(\NvN,\Gamma,\beta)$ giving
rise to $\MvN\simeq\NvN\rtimes_\beta\Gamma$ be the discrete decomposition
of $\MvN$ associated to $\phi$.
Let $(\NvNt,\Rp,\betat)$ giving rise to $\MvN\simeq\NvNt\rtimes_\betat\Rp$
be Takesaki's continuous decomposition.
Then the action $\betat$ of $\Rp$ on $\NvNt$ is the action induced
up to $\Rp$ from the action $\beta$ of $\Gamma$ on $\NvN$,
(where $\Gamma$ is embedded in $\Rp$ by $\ioh$).
\endproclaim
\demo{Proof}
Apply Proposition 3.3 with $\alpha$ and $\iota$ as in \S1.3.
\QED

\noindent{\bf \S4. A model for the one--parameter trace--scaling action
on the injective II$_\infty$ factor.}

\proclaim{Lemma 4.1}
Let $\NvN_i$ be a von Neumann algebra with almost periodic
n.f.s\. weight $\psi_i$, ($i=1,2$).
Let $\MvN=\NvN_1\otimes\NvN_2$ be the tensor product von Neumann algebra
and $\phi=\psi_1\otimes\psi_2$ the tensor product weight.
Then $\phi$ is an almost periodic n.f.s. weight on $\MvN$ and
the point spectrum of $\Delta_\phi$ is the product of the point spectra
of $\Delta_{\psi_1}$ and $\Delta_{\psi_2}$.
Let $\Gamma_i$ (respectively $\Gamma$), taken with discrete topology,
be the subgroup of $\Rp$ generated by the point spectrum of
$\Delta_{\psi_i}$ (respectively $\Delta_\phi$),
so $\Gamma=\Gamma_1\Gamma_2$, and let $G_i=\widehat\Gamma_i$
(respectively $G=\widehat\Gamma$).
Let $\alpha_i$ (respectively $\alpha$) be the action of $G_i$ on $\NvN_i$
(respectively $G$ on $\MvN$) as in \S1.3.
Then for $\gamma\in\Gamma$, the spectral subspace is
$$ \MvN^\alpha(\{\gamma\})
=\bigoplus_{\{\gamma_1\in\Gamma_1,\gamma_2\in\Gamma_2
\mid\gamma_1\gamma_2=\gamma\}}
\NvN_1^{\alpha_1}(\{\gamma_1\})\otimes\NvN_2^{\alpha_2}(\{\gamma_2\}), \tag9 $$
meaning the weak closure of the set of linear
combinations of simple tensors $a_1\otimes a_2$ with
$a_i\in\NvN_i^{\alpha_i}(\{\gamma_i\})$, ($i=1,2$) where
$\gamma_1\gamma_2=\gamma$.
\endproclaim
\demo{Proof}
{}From $\L^2(\MvN,\phi)=L^2(\NvN_1,\psi_1)\otimes L^2(\NvN_2,\psi_2)$
and $\Delta_\phi=\Delta_{\psi_1}\otimes\Delta_{\psi_2}$ one sees
that $\phi$ is almost periodic and has the desired point spectrum.
Moreover, from the characterization in Lemma~1.3.3, one sees that
$\supseteq$ holds in~(9).
But by Lemma~1.2.3, the direct sum of $\NvN_1^{\alpha_1}(\{\gamma_1\})
\otimes\NvN_2(\{\gamma_2\})$ over all $\gamma_1\in\Gamma_1$ and
$\gamma_2\in\Gamma_2$ is dense in $\MvN$, so equality must hold
in~(9).
\QED

\proclaim{Theorem 4.2}
Let $0<\mu_1,\mu_2<1$ be such that $\mu_1/\mu_2$ is irrational.
Let $\Gamma$ be the multiplicative
subgroup of $\Rp$ generated by $\mu_1$ and $\mu_2$,
let $\NvN$ be the injective II$_\infty$ factor and let $\beta$
be an action of $\Gamma$ on $\NvN$ such that $\Tr_\NvN\circ\beta_\gamma
=\gamma\Tr_\NvN$ $\forall\gamma\in\Gamma$, where $\Tr_\NvN$ is the
n.f.s\. trace on $\NvN$.
Then the action of $\Rp$ induced up from the action $\beta$ of $\Gamma$
on $\NvN$ is the one--parameter trace--scaling action of $\Rp$
on the injective II$_\infty$ factor.
\endproclaim

\proclaim{Remark 4.3}\rm
Since a trace--scaling action $\beta$ of $\Gamma$ on $\NvN$ is relatively
easy to construct, Theorem~4.2 gives an accessible picture
of the one--parameter trace--scaling action on the
injective II$_\infty$ factor, heuristically
as translation on the fixed point subalgebra of $L^\infty(\Rp,\NvN)$
under a pair of commuting automorphisms (see Definition 3.1)
\endproclaim

\demo{Proof of Theorem 4.2}
By~\cite{11}, all trace--scaling actions $\beta$ of $\Gamma$ on $\NvN$
are outer conjugate, so it suffices to show the theorem for a
particular one.
Let $\lambda_i=(1+\mu_i)^{-1}$ and let $\MvN^{\lambda_i}$ be Powers'
factor with Powers' state $\omega_i=\omega_{\lambda_i}$ \cite{12}.
It is well--known ({\it cf\.} \cite{13}, \S5) that the centralizer
$\MvN^{\lambda_i}_{\omega_i}$ is the hyperfinite II$_1$--factor, $R$,
and the point spectrum of $\Delta_{\omega_i}$ equals
$\{\mu_i^n\mid n\in\Integers\}$.
Consider the AFD (approximately finite dimensional), purely infinite factor
$\MvN=\MvN^{\lambda_1}\otimes\MvN^{\lambda_2}$ and let
$\phi=\omega_1\otimes\omega_2$ be the tensor product state.
Then by Lemma~4.1
$\phi$ is almost periodic, the point
spectrum of $\Delta_\phi$ is $\Gamma=\{\mu_1^{n_1}\mu_2^{n_2}\mid
n_1,n_2\in\Integers\}$ and the centralizer
$\MvN_\phi=R\otimes R\simeq R$ is a factor.
Now by Proposition 2.12, the discrete decomposition of $\MvN$ associated
to $\phi$ is a trace--scaling action of $\Gamma$ on the injective
II$_\infty$ factor.
By the factoriality of the centralizer and
Connes~\cite{4}, 3.2.7\., it follows that
$\MvN$ is a type III$_1$ factor, so Takesaki's continuous
decomposition is the one--parameter trace--scaling action of $\Rp$
on the injective type II$_\infty$ factor.
Now Corollary~3.4 finishes the proof.
\QED

\proclaim{Remark 4.4}\rm
A more concrete picture of the one--parameter trace--scaling
action in terms of operators on Hilbert space can be had from
the proof of Theorem~4.2.
In~\cite{13}, \S5, Takesaki showed that $\omega_i$ is a homogeneous,
periodic state.
Thus, by general results of the same paper,
there is an isometry, $v_i$, in the spectral subspace
$\MvN^{\lambda_i}_{\omega_i}(\{\mu_i^{-1}\})$ such that
$$ \MvN^{\lambda_i}_{\omega_i}(\{\mu_i^n\})
=\cases M^{\lambda_i}_{\omega_i}v_i^{-n}&\text{ if }n\le0 \\
(v_i^*)^nM^{\lambda_i}_{\omega_i}&\text{ if }n\ge0. \endcases \tag10 $$
Keeping in mind that the centralizers $\MvN^{\lambda_i}_{\omega_i}$
are copies of the hyperfinite II$_1$ factor, $R$, and noting that
the spectral subspaces of the tensor product are
$$ \MvN_\phi(\{\mu_1^{n_1}\mu_2^{n_2}\})
=\MvN^{\lambda_1}_{\omega_1}(\{\mu_1^{n_1}\})\otimes
\MvN^{\lambda_2}_{\omega_2}(\{\mu_2^{n_2}\}), $$
one can use~(10) to describe the spectral subspaces of $\MvN$
under $\phi$.
Let $\Hil$ be the Hilbert space on which $\MvN$ acts.
Then by Proposition~3.2, the injective II$_\infty$ factor
is the von Neumann algebra on $L^2(\Rp)\otimes\Hil$ generated by
$$ (L^\infty(\Rp)\otimes1)\cup\{\lambda_\gamma\otimes a
\mid \gamma\in\Gamma,\,a\in\MvN_\phi(\{\gamma\})\} $$
and the one--parameter trace--scaling action is given by
$\Ad(\lambda_t\otimes1)$ for $t\in\Rp$, {\it i.e\.} translation
by $t$ on the first component.
\endproclaim

\noindent{\bf\S5. Another model for the one--parameter
trace--scaling action on the injective II$_\infty$ factor.}

  We now explain a model which was described to the author by
V.F.R\. Jones, and depends on work of P.-L\. Aubert.
The usual (linear) action of
$SL(2,\Integers)$ on $\Real^2$ preserves Lebesque measure.
Let $\NvN$ be the crossed product von Neumann algebra
of $L^\infty(\Real^2)$ by the  action of $SL(2,\Integers)$.
{}From Lebesque measure on $\Real^2$, one gets a n.f.s\.
trace on $\NvN$.
In~\cite{2}, Aubert showed the somewhat surprising fact that $\NvN$
is the injective II$_\infty$ factor.
Then for $t\in\Rp$, the dilation of $\Real^2$
by $t$ gives an automorphism, $\alpha_t$,
of $\NvN$ that scales the trace by $t$.
So $t\mapsto\alpha_t$ is the one--parameter trace--scaling
action of $\Rp$ on the injective II$_\infty$ factor.

\enddocument